\documentclass[12pt]{amsart}
\usepackage{amssymb, amsthm, pb-diagram, pstricks, pst-tree}

\newtheorem{theorem}{Theorem}[section]
\newtheorem{proposition}[theorem]{Proposition}
\newtheorem{thm}{Theorem}[subsection]
\newtheorem{prop}[thm]{Proposition}
\newtheorem{cor}[thm]{Corollary}
\theoremstyle{definition}
\newtheorem{defn}[thm]{Definition}
\newtheorem{rem}[thm]{Remark}
\newtheorem{remark}[theorem]{Remark}
\numberwithin{equation}{section}

\title[Infinity Structure of Poincar\'e Duality Spaces]{Infinity Structure of Poincar\'e Duality Spaces}

\author[T.Tradler \and M. Zeinalian]{Thomas Tradler and Mahmoud Zeinalian}

\address{Dennis Sullivan, Department of Mathematics,
Graduate Center of the City University of New York, 365 Fifth Avenue, New
York, NY 10016, USA.} \email{dsullivan@gc.cuny.edu}

\address{Thomas Tradler, Department of Mathematics, College of
Technology of the City University of New York, 300 Jay Street,
Brooklyn, NY 11201, USA.} \email{ttradler@citytech.cuny.edu}

\address{Mahmoud Zeinalian, Department of Mathematics, C.W. Post Campus
of Long Island University, 720 Northern Boulevard, Brookville, NY
11548, USA.} \email{mzeinalian@liu.edu}
\begin{document}

\maketitle \vspace{-0.6cm} \centerline{with an appendix by {\sc
Dennis Sullivan}}

\begin{abstract}
We show that the complex $C_\bullet X$ of rational simplicial chains
on a compact and triangulated Poincar\'e duality space $X$ of
dimension $d$ is an A$_\infty$ coalgebra with $\infty$ duality. This
is the structure required for an A$_\infty$ version of the cyclic
Deligne conjecture. One corollary is that the shifted Hochschild
cohomology $HH^{\bullet+d} (C^\bullet X, C_\bullet X)$ of the
cochain algebra $C^\bullet X$ with values in $C_\bullet X$ has a BV
structure. This implies, if $X$ is moreover simply connected, that
the shifted homology $H_{\bullet+d}LX$ of the free loop space admits
a BV structure. An appendix by Dennis Sullivan gives a general local
construction of $\infty$ structures.
\end{abstract}


\section{Introduction}
\subsection{Introduction} In \cite{TZ1}, it was shown that the Hochschild complex of an
associative algebra with an invariant, symmetric, and
nondegenerate inner product has a natural action of the PROP of
Sullivan chord diagrams. This was referred to as the cyclic
Deligne conjecture. Unfortunately, such a statement does not
directly apply to the topology of Poincar\'e duality spaces, since
the associative algebra of simplicial cochains does not enjoy a
nondegenerate and symmetric inner product which is also invariant.
As it turns out, the right notion for this cyclic Deligne
conjecture is that of an A$_\infty$ coalgebra with $\infty$
duality; that is, an A$_\infty$ coalgebra $C$ together with a
duality map $F$ enjoying appropriate invariance, nondegeneracy,
and symmetry properties; see Definition \ref{duality}. This
structure was defined in \cite{T2}, where it was shown that the
Hochschild cohomology $HH^{\bullet+d}(C^\ast, C)$ is a BV algebra.
In fact, it was later shown in \cite{TZ2} that the PROP
$\mathcal{DG}_2$, an enhancement of the PROP of Sullivan chord
diagrams, acts on this Hochschild complex $CH^{\bullet+d} (C^\ast
,C)$. In this paper, we give an explicit construction of an
$\infty$ duality structure on $C_\bullet X$, the rational chains
on a compact and triangulated Poincar\'e duality space $X$ of
dimension $d$. This implies the shifted Hochschild cohomology
$HH^{\bullet+d} (C^\bullet X, C_\bullet X)$ of the cochains
$C^\bullet X$ with values in $C_\bullet X$, has a natural BV
structure with a unit. A corollary of this, when $X$ is in
addition simply connected, is that the shifted homology
$H_{\bullet+d}LX$ of the free loop space admits a unital BV
structure. This should be viewed in light of Chas and Sullivan's
seminal work on the algebraic topology of the free loop space
$LX$, when $X$ is a manifold \cite{CS1}.

Let us recall some of the history and motivations. Chas and
Sullivan \cite{CS1} defined a BV structure on the shifted homology
$H_{\bullet+d}LX$ of the free loop space $LX$ of an orientable
manifold $X$, with a unit only in the closed case; see Remark
\ref{unit}. More precisely, they defined a multiplication,
$$\centerdot: H_{\bullet+d} LX \otimes H_{\bullet+d} LX \to
H_{\bullet+d}LX,$$ called the loop product, along with an
operator, $$\Delta : H_{\bullet+d}LX \to H_{\bullet+1+d}LX,$$ and
proved that the triple $(H_{\bullet+d}LX , \centerdot, \Delta)$
forms a BV algebra. That is to say that $(H_{\bullet+d} LX,
\centerdot)$ is an associative and graded commutative algebra with
a unit along with a differential $\Delta$ whose deviation form
being a derivation defines a structure of a graded Lie algebra on
$H_{\bullet+d+1}LX$. The multiplication, $\centerdot$, owes itself
to the fact that two loops with common base points can be composed
and that in a manifold cycles can be made transversal and then
intersected. They coined the term {\it string topology} to
describe the activity aimed towards understanding the algebraic
structure of the free loop space of a manifold. Since then, there
have been several generalizations and applications of their
machinery in different directions; see, for instance, Sullivan
\cite{DS2}, Chas \& Sullivan  \cite{CS2}, Cohen \& Klein \&
Sullivan \cite{CKS}, Cohen \& Jones \cite{CJ}, Klein
\cite{Klein2003}, Felix \& Thomas \& Vigue-Poirrier \cite{FTV},
Merkulov \cite{Me}, and Chataur \cite{Chataur2003}.

\begin{rem} \label{Klein-PD}
Extension of string topology to more general spaces than manifolds
have been considered by others. For example, J. Klein
\cite{Klein2003} defined a product on the loop space homology
$H_\bullet LX$ for an orientable Poincar\'e duality space $X$. In
the case of a manifold, this product coincides with that of String
Topology.
\end{rem}

\begin{rem}\label{HomotopyInv}
It has been shown in \cite{CKS} that the Chas-Sullivan BV
structure is in fact a homotopy invariant of the manifold. This fits
with the work of this paper as well as with the view that the
Chas-Sullivan BV structure is defined for a lager category of spaces
than manifolds. Nonetheless, Sullivan conjectures that the further
string topology, involving operations labelled by rest of moduli
space of Riemann surfaces, is not a homotopy invariant; see
Conjecture 3 in the Postscript of \cite{DS2}.
\end{rem}

A corollary of the result of this paper is that there is a unital
BV structure on the shifted cohomology of the larger class of
simply connected, compact, and triangulated Poincar\'e duality
spaces. By definition, a Poincar\'e duality space $X$ of dimension
$d$ is a topological space with an element $[X] \in H_dX$, called
the fundamental class, such that the map $[X] \cap \cdot :
H^\bullet X \to H_{d-\bullet} X$ is an isomorphism. A subtle but
very important point is that although the above map is a map of
$H^\bullet X$ bimodules, when viewed at the chain level as a map
$C^\bullet X \to C_{d-\bullet} X$, is not a map of $C^\bullet
X$-bimodules. A closed and oriented manifold is a good example of
a Poincar\'e duality space and it is therefore natural to ask if
in this case our structure coincides with the Chas-Sullivan BV
structure. This question, however, has not been addressed in this
paper.

Let us highlight some of the salient points of this paper. We have
used an inductive and combinatorial argument to construct a
symmetric A$_\infty$ co-inner product, and subsequently used
minimal models to obtain an $\infty$ duality structure. As a
result we show,

\begin{itemize}

\item[i)]{The chains $C_\bullet X$ on a compact and triangulated
Poincar\'e duality space $X$ is naturally an A$_\infty$ coalgebra
with $\infty$ duality.}

\item[ii)]{For an $X$ as above, $HH^{\bullet+d}(C^\bullet X,
C_\bullet X)$ is a BV algebra with a unit; see Remark \ref{unit}.}

\item[iii)]{If $X$ is in addition simply connected, then the
shifted homology $H_{\bullet+d}LX$ of the free loop space is
naturally a BV algebra with a unit.}

\end{itemize}

\begin{rem} \label{unit}
Whether or not the Chas-Sullivan loop multiplication has a unit is
decided by whether or not the manifold in question is closed, i.e.
compact and without boundary. This suggest there should be a
version of the construction of this paper, without a unit, which
works for the larger category of Lefschetz duality spaces.
\end{rem}

A general and very useful scheme for the local construction of
$\infty$ structures is provided by Dennis Sullivan in Appendix
\ref{local-sonstruction-DS}. In Appendix \ref{circ}, an explicit
calculation of the A$_\infty$ co-inner product for the case of the
circle shows the necessity of having higher homotopies even when
the A$_\infty$ coalgebra structure is in fact strictly
coassociative. Such an explicit calculation of the structure in
higher dimensions, even for the case of a two dimensional compact
and orientable surface, should be very interesting.

Let us review the content of each chapter in more details. In
Chapter 2, we present the algebraic definitions necessary for the
future chapters. An A$_\infty$ coalgebra structure over a graded
module $C$ is a differential on the free associative algebra $BC$
over the shifted vector space $C$; see \cite{S}. An A$_\infty$
cobimodule over an A$_\infty$ coalgebra is defined in a similar
fashion by having a differential satisfying certain properties in
terms of its interaction with the differential of the A$_\infty$
coalgebra over which it is defined; see \cite{M}. An $\infty$
duality is a map $F$ between A$_\infty$ cobimodule $C$ and its
dual $C^\ast$, which satisfies a symmetry condition and induces an
isomorphism at the homology level; see \cite{TT}. We recall from
\cite{T2} that an $\infty$ duality structure on $C$ yields a BV
structure on the Hochschild cohomology $HH^\bullet (C^*, C)$. In
Chapter 3, we show that the chains $C_\bullet X$ on a compact
polyhedron $X$ of dimension $d$, together with a choice of a cycle
$\mu \in C_\bullet X$, is naturally furnished with a locally
constructed A$_\infty$ co-inner product structure satisfying a
symmetry condition; see Theorem \ref{IPonManifold}. Moreover, we
show that when $X$ is a Poincar\'e duality space and $\mu=[X]$ its
fundamental class, this symmetric A$_\infty$ co-inner product is
an $\infty$ duality structure; that is, the A$_\infty$ co-inner
product structure induces an isomorphism at the level of the
homology; see Theorem \ref{non-deg}. This follows from Proposition
\ref{invertability} which is proved in two steps. First, we argue
that an A$_\infty$ module over an A$_\infty$ algebra has a minimal
model; see subsection \ref{minimal}. Then, using method similar to
those employed in \cite{DS1} and \cite{FHT}, we show that an
A$_\infty$ module map between two minimal models which induces an
isomorphism at the level of the homology is in fact an
isomorphism. Algebraically, this inverse map can be written down
very explicitly. It is noteworthy that even though the $\infty$
duality map is local, there is no guarantee for the locality of
the algebraically constructed inverse map; see Remark
\ref{LocInv}. Nonetheless, such a structure induces a
quasi-isomorphism between the Hochschild complexes $CH^\bullet
(C^\bullet X, C^\bullet X)$ and $CH^\bullet (C^\bullet X,
C_\bullet X)$. According to \cite{T2} the transport of the
algebraic structures from one Hochschild cohomology onto the
other, and putting everything together, after a shift in degrees,
yields a BV structure. A result of Jones \cite{J} states that, for
a simply connected space $X$, the Hochschild cohomology
$HH^\bullet (C^\bullet X, C_\bullet X)$ calculates the homology
$H_\bullet LX$ of the free loop space. Therefore, it follows that
if the Poincar\'e duality space $X$ is in addition simply
connected, the shifted homology $H_{\bullet+d} LX$ of the free
loop space $LX$ is naturally a BV algebra. There are two
appendices to this paper. Appendix \ref{local-sonstruction-DS},
written by Dennis Sullivan, provides a general local construction
for infinity structures. In a second appendix, we study the
$\infty$ duality structure on $C_\bullet S^1$, the chains on the
circle $S^1$. This structure, in terms of its component maps, is
explicitly calculated.

\textbf{Acknowledgements} We would like to thank Dennis Sullivan
for adding an appendix to this paper and for his continuous
support and suggestions. We would like to thank Andrew Ranicki and
Ralph Cohen for useful discussions, and Scott Wilson for pointing
out a gap in an earlier version. We are also thankful for the
referee's comments which have improved the presentation of the
paper.
\section{Algebraic definitions}

\subsection{$\infty$ duality structure}

In this section, we recall some of the pertinent algebraic and
homotopy notions. The Poincar\'e duality constitutes a $H^\bullet
X$-bimodule isomorphism between cohomology $H^\bullet X$ and
homology $H_{d-\bullet}X$ of a Poincar\'e duality space $X$ of
dimension $d$. The correct homotopy notion is obtained by resolving
the operad of bimodule maps.

\begin{rem}
It has been shown in \cite{LT} that given a cyclic quadratic
operad $\mathcal O$, there exists a colored operad $\hat{\mathcal
O}$ so that an algebra over $\hat{\mathcal O}$ is precisely an
algebra over $\mathcal O$ together with a symmetric and invariant
inner product. Furthermore, the Koszulness of $\mathcal O$ as an
operad implies Koszulness of $\hat{\mathcal O}$ as a colored
operad. More precisely, let $\mathcal O^!$ denote the quadratic
dual operad of $\mathcal O$ and $\textbf{D}(\mathcal O^!)$ its
cobar dual operad. Also, let $\hat{\mathcal O}^!$ be the quadratic
dual colored operad of $\hat{\mathcal O}$, whose cobar dual
colored operad is denoted by $\textbf{D}(\hat{\mathcal O}^!)$. The
Koszulness of $\mathcal O$ means that the canonical map
$\textbf{D}(\mathcal O^!)\to \mathcal O$ is a quasi-isomorphism of
operads. This implies that the canonical map $\textbf{D}( \hat{
\mathcal O}^!)\to \hat{ \mathcal O}$ is also a quasi-isomorphism
of colored operads; see Theorem 2.8 of \cite{LT}. This statement
allows one to resolve $\hat{\mathcal O}$ to obtain a homotopy
version $\hat{ \mathcal O}_\infty= \textbf{D}( \hat{ \mathcal
O}^!)$ of both the algebra and the inner product.
\end{rem}

We now apply the result from \cite{LT} to the cyclic quadratic
Koszul operad $\mathcal Assoc$ which governs associative algebras.
Let us describe the outcome of the construction $\widehat{\mathcal
Assoc}_\infty$ more explicitly. Let $C=\bigoplus_{j\in \mathbb{Z}}
C_{j}$ be a graded module over fixed associative and commutative
ring $R$ with a unit. Define the suspension $sC$ of $C$ as the
graded module $sC= \bigoplus_{j\in \mathbb{Z}} (sC)_{j}$ with
$(sC)_{j}= C_{j+1}$. The suspension map $s:C\to sC$, $v\mapsto
sc=c$ is a linear isomorphism of degree $-1$. We denote by
$BC=T(sC)=\prod_{i\geq 0} (sC)^{\otimes i}$ the (completed) free
associative algebra on the suspended space $sC$.

An A$_\infty$ coalgebra structure on $C$ is defined to be a
derivation $D\in Der(BC)$ of degree $+1$ with $D^{2}=0$. The space
$Der(BC)$ has a differential $\delta:Der(BC) \to Der(BC)$ given by
$\delta (D')=[D,D']=D\circ D'-(-1)^{|D'|} D' \circ D$ with
$\delta^{2}= 0$. Given an A$_\infty$ coalgebra $(BC,D)$ and a
graded $R$-module $M$, let $B^MC=T^{sM}(sC)= \prod_{i,j\geq 0} (sC
)^{\otimes i}\otimes (sM)\otimes (sC)^{\otimes j}$. An A$_\infty$
cobimodule structure on $M$ over $C$ is by definition a derivation
$D^{M}\in Der(B^MC)$ over $D$ of degree $+1$ with $(D^{M})^{2}=0$;
see \cite{TT}.

For A$_\infty$ cobimodules $(B^M C,D^{M})$ and $(B^N C,D^N)$ over
the A$_\infty$ coalgebra $BC$, let $Hom(B^MC, B^NC)$ denote the
space of all module maps. $Hom(B^MC, B^NC)$ has a natural
differential given by $\delta^{M,N}(F)=D^N\circ F-(-1)^{|F|}F\circ
D^{M}$. An A$_\infty$ cobimodule map between $B^M C$ and $B^N C$
is by definition an element $F\in Hom(B^MC,B^NC)$ of degree 0 with
$\delta^{M,N}(F)=0$; that is, $ D^N\circ F=F\circ D^{M}$. Note
that since $F$ is a map between free objects, it is uniquely
determined by the components $F|_M:M\to B^N C$.

\begin{rem} The concepts of A$_\infty$ coalgebras, A$_\infty$
cobimodules, and A$_\infty$ cobimodule maps are generalizations of
those of coassociative coalgebras, cobimodules, and cobimodule
maps; see \cite{TT}.
\end{rem}

Consider an A$_\infty$ cobimodule map $F:B^M C\to B^N C$. We will
argue that the induced map $F_*:H_\bullet(B^M C,D^M)\to
H_\bullet(B^N C,D^N)$, between the homologies, is an isomorphism
if and only if the lowest component $F_{0,0}:M\to N$ of $F$
induces an isomorphism $(F_{0,0})_*: H_\bullet(M, D^M_{ 0,0})\to
H_\bullet(N,D^N_{0,0})$ between the homologies. Such a map is
called a \emph{quasi-isomorphism}. We show this using by using a
minimal model decomposition for A$_\infty$ modules.

\subsection{Minimal Models and Invertibility}\label{minimal}

It was shown in \cite{DS1} that every free connected differential
graded algebra is isomorphic to the tensor product of a unique
minimal algebra and a unique contractible algebra.  Since, similar
decomposition theorems have been proved and utilized in different
contexts; for instance see section 8 of \cite{FHT}, Theorem 5.2 of
\cite{Ka}, and the Lemma in section 4.5.1 of \cite{Ko}. In this
section, we prove that every A$_\infty$ cobimodule $(B^M C,D^M)$
over $(BC,D)$ can be decomposed into a direct sum of a linear
contractible and a minimal one. This will allow us to show that an
A$_\infty$ cobimodule map $F:B^M A\to B^N A$ induces an isomorphism
on homology $F_*:H_\bullet(B^M A,D^M)\to H_\bullet(B^N A,D^N)$, if
the lowest component $F_{0,0}:M\to N$ induces an isomorphism on
homology $(F_{0,0})_*: H_\bullet(M,D^M_{0,0})\to H_\bullet(N,D^N_{0,
0})$. For similar a statement in a different context see Theorem 5.4
of \cite{Ka}. For a more general operadic method see (M3) of \cite{M2}.

\begin{defn}\label{quasi-iso} Let $(BC,D)$ be an A$_\infty$ coalgebra
and $(B^P C,D^P)$ be an A$_\infty$ cobimodule. $(B^P C,D^P)$ is
called \emph{minimal}, if the lowest component $D^P_{0,0}:M\to M$
of the map $D^P:B^P C\to B^P C$  vanishes, i.e. $D^P_{0,0}=0$.
\end{defn}

\begin{defn} Let $(BC,D)$ be an A$_\infty$ coalgebra
and $(B^N C,D^N)$ be an A$_\infty$ cobimodule. $(B^N C,D^N)$ is
called \emph{linear contractible}, if the lowest component of
$D^N:B^N C\to B^N C$ is the only non-vanishing component, i.e. for
$i+j>0$ it is $D^N_{i,j}=0$, and additionally the homology of
$(N,D^N_{0,0})$ vanishes.
\end{defn}

With this we can state the decomposition theorem.
\begin{thm}\label{decomposition} Let $(BC,D)$ be an $A_\infty$
coalgebra. Every $A_\infty$ cobimodule $(B^M C, D^M)$ is
isomorphic as an A$_\infty$ cobimodule to the direct sum of a
minimal one $(B^P C,D^P)$ and a linear contractible one $(B^N
C,D^N)$. That is, $$ B^M C \cong B^P C \oplus B^N C.$$
\end{thm}
\begin{proof}
Choose an arbitrary positive definite inner-product on the vector
space $M$. Denote by $d$ the lowest component of the differential
on $B^M C$, i.e. $d=D^M_{0,0}$. Define $d ^\dagger$ to be the
adjoint of $d$ under the chosen inner product. Then, define
$X=Im(d)$, $Y=Im(d^\dagger)$, and the harmonic subspace
$P=Ker(d)\cap Ker(d^\dagger)$. One knows that $M$ has a Hodge
decomposition $M=P\oplus X\oplus Y$. We will use the decomposition
$M=P\oplus N$, where $N=X\oplus Y$. An elementary check shows that
$d:Y\to X$ is an isomorphism of vector spaces and that every
homology class $[m]$ of $(M,d)$ has a unique harmonic
representative.

We want to show that there exists a sequence of $A_\infty$
cobimodules structures $\{ (B^{P\oplus N} C,D^{P\oplus N}(n))
\}_{n\in \mathbb N _0}$ on the space $B^{P\oplus N} C$, whose
components will be written as $D^{P\oplus N}(n)|_{P\oplus
N}=d^n_0+d^n_1 +d^n_2+\cdots$, where $ d^n_j:P\oplus N\to
(B^{P\oplus N} C)_j= \bigoplus_{r+s=j} C^{\otimes r}\otimes
(P\oplus N)\otimes C^{\otimes s}$. These maps will satisfy,
\begin{eqnarray}\label{induction0}
\quad\quad d^n_0(x)=0 &\quad d^n_0(y)\in X &\quad d^n_0(p)=0,  \\
 \label{inductionk}
\quad\quad d^n_k(x)=0 &\quad d^n_k(y)=0    &\quad d^n_k(p)\in (B^P
C)_k,
\end{eqnarray}
for $x\in X$, $y\in Y$, $p\in P$, and $1\leq k\leq n$, together with the compatibility
condition $$ d^n_k=d^m_k,  \quad \text{for} \quad k\leq min(m,n).$$

Furthermore, there are $A_\infty$ cobimodules isomorphisms,
$$\varphi(n):(B^M C,D^M) \to (B^{P\oplus N} C,D^{P\oplus N} (n)),$$
satisfying the compatibility,$$ \varphi^n_k=\varphi^m_k, \quad\text{for}\quad
k\leq min(m,n),$$ where $\varphi^n_k$ is the component of
$\varphi(n)|_M =\varphi^n_0+\varphi^n_1+\varphi^n_2+ \cdots$ that
maps to $B^{(P\oplus N)} C_k$. After constructing these maps and
structures, we obtain an $A_\infty$ cobimodule structure on
$B^{(P\oplus N)} C$ by taking components $(D^{P\oplus N})_k=d^n_k$
for any $n\geq k$. Using equation (\ref{induction0}) and
(\ref{inductionk}), it is clear that this splits into a minimal
$A_\infty$ cobimodule $(B^P C, D^P)$ and a linear contractible one
$(B^N C,D^N)$. Furthermore the maps $\varphi(n)$ induce an
$A_\infty$ cobimodule isomorphism $\varphi:B^M C\to B^P C\oplus
B^N C=B^{(P\oplus N)} C$ by $\varphi_k=\varphi^n_k$, for $n\geq
k$.

Let us start with the case $n=0$. We identify $B^{P\oplus N} C=B^M
C$, and take $D^{P\oplus N}(0)=D^M$ together with $\varphi(0)=id:
B^M C\to B^{P\oplus N} C$. Then, condition (\ref{induction0}) is
satisfied, because we have $d^0_0=d$, and with the definitions of
$X$, $Y$ and $P$ it is clear that $d(X)=d(P)=\{0\}$, while
$d(Y)\subset X$.

Now, let us assume we have constructed $D^{P\oplus N}(k)$ and
$\varphi(k)$ for $k=1,\cdots,n$, satisfying (\ref{induction0}) and
(\ref{inductionk}). Then, define a tensor-bimodule isomorphism
$F:B^{P\oplus N} C\to B^{P\oplus N} C$ to be given by $F|_{P\oplus
N }=id+f$. The map $f:P\oplus N\to (B^{P\oplus N} C)_{n+1}$ is
given by,
\begin{eqnarray*}
 &\text{for} \quad x\in X& \quad f(x)=-d^n_{n+1}((d^n_0)^{-1}(x)),\\
 &\text{for} \quad y\in Y& \quad f(y)=0,\\
 &\text{for} \quad p\in P& \quad f(p)\quad \text{is so that}\quad
       \widetilde{d^n_0}(f(p))=(id-\pi)\circ d^n_{n+1}(p),
\end{eqnarray*}
where $d^n_0=d:Y\to X$, and $\pi$ is the projection $\pi:B^M C\to
B^P C$. In order to see that $f(p)$ can be defined, first notice
that $0=(D^{P\oplus N}(n))^2(p)=\widetilde{d^n_0}\circ
d^n_{n+1}(p)+ \sum_{j=1}^{n+1}\widetilde{ d^n_j}\circ
d^n_{n+1-j}(p)$. By the inductive hypothesis (\ref{induction0})
and (\ref{inductionk}) it is $\widetilde{d^n_0}\circ
d^n_{n+1}(p)\in B^N C$ and $\sum_{j=1}^{n+1}
\widetilde{d^n_j}\circ d^n_{n+1-j}(p)\in B^P C$, so that
$\widetilde{d^n_0}\circ d^n_{n+1}(p)=0$. Furthermore, it is also
$\widetilde{d^n_0}\circ \pi \circ d^n_{n+1}(p)=\pi\circ\widetilde{
d^n_0}\circ d^n_{n+1}(p)= 0$. Thus, we see that $(id-\pi)\circ
d^n_{n+1}(p)\in B^N C$ is a closed element. But since
$H_\bullet(B^N C,\widetilde{d^n_0})=\cdots\otimes
H_\bullet(N,d|_N)\otimes\cdots=\{0\}$, one can find an element
$f(p)\in B^N C$ so that $\widetilde{ d^n_0} (f(p))=
(id-\pi)\circ d^n_{n+1}(p)$.

Thus, $F$ defines a vector space isomorphism with inverse given by
$F^{-1}|_{P\oplus N}=id-f+ \text{higher order terms}$. With this,
we define $D^{P\oplus N}(n+1)$ to be the induced $A_\infty$
cobimodule structure of $F$ on $B^{P\oplus N} C$, i.e. $D^{P\oplus
N}(n+1)=F\circ D^{P\oplus N}(n)\circ F^{-1}$. This defines an
$A_\infty$ cobimodule structure and we get an $A_\infty$
cobimodule isomorphism by taking $\varphi(n+1)=F\circ \varphi(n)$.
Notice that $D^{P\oplus N}(n+1)=F\circ D^{P\oplus N}(n)\circ
F^{-1}=(id+f)\circ (d^n_0+d^n_1+\cdots)\circ(id-f+\cdots)$, so
that $d^{n+1}_k=d^n_k$ for $k\leq n$, and
$d^{n+1}_{n+1}=d^n_{n+1}+ f\circ d^n_0- d^n_0 \circ f$. Similarly
$\varphi(n+1)=F \circ \varphi(n)= (id+f)
\circ(\varphi^n_0+\varphi^n_1+ \cdots)$ implies $\varphi
^{n+1}_k=\varphi^n_k$ for $k\leq n$. Since the differential
$D^{P\oplus N}(n)$ and the isomorphism $\varphi(n)$ are not
altered up to the $n^{th}$ level, only equation (\ref{inductionk})
with $k=n+1$ requires a check. Recall that $d^n_0=d^0_0=d$, so
that for $x\in X=Im(d)$, we have,
\begin{eqnarray*} \quad\quad\quad
 d^{n+1}_{n+1}(x)&=&d^n_{n+1}(x)+f(d^n_0(x))-\widetilde{d^n_0}(f(x))\\
   &=&d^n_{n+1}(x)+\widetilde{d^n_0}(d^n_{n+1}((d^n_0)^{-1}(x)))\\
   &=&d^n_{n+1}(x)-d^n_{n+1}(d^n_0((d^n_0)^{-1}(x)))\\
   &=&0.
\end{eqnarray*}
Note that we made use of the fact that $\left(D^{P\oplus
N}(n)\right)^2(y)=0$, for $y\in Y$, implies $0= d^n_{n+1}\circ
d^n_0(y)+\widetilde{ d^n_{n}}\circ d^n_1(y)+\cdots+
\widetilde{d^n_0} \circ d^n_{n+1}(y)= d^n_{n+1}\circ d^n_0(y)+
\widetilde{d^n_0}\circ d^n_{n+1}(y)$, because $d^n_1(y)=
d^n_2(y)=\cdots=d^n_n(y)=0$ by condition (\ref{inductionk}).

Next, for any $y\in Y$, we have,
\begin{eqnarray*} \quad\quad\quad
 d^{n+1}_{n+1}(y)&=&d^n_{n+1}(y)+f(d^n_0(y))-\widetilde{d^n_0}(f(y))\\
   &=&d^n_{n+1}(y)-d^n_{n+1}((d^n_0)^{-1}(d^n_0(y)))\\
   &=&0,
\end{eqnarray*}
because $d^n_0(y)\in X$.

Finally, if $p\in P$, we have,
\begin{eqnarray*} \quad\quad\quad
 d^{n+1}_{n+1}(p)&=&d^n_{n+1}(p)+f(d^n_0(p))-\widetilde{d^n_0}(f(p))\\
   &=&d^n_{n+1}(p)-(id-\pi)\circ d^n_{n+1}(p)\\
   &=&\pi\circ d^n_{n+1}(p)\in B^P C.
\end{eqnarray*}
\end{proof}

\begin{defn} Let $(BC,D)$ be an A$_\infty$ coalgebra
and let $F:(B^{M_1} C,D^{M_1})\to (B^{M_2} C,D^{M_2})$ be an
A$_\infty$ cobimodule morphism between two A$_\infty$ cobimodules.
Then, $F$ is called a \emph{quasi-isomorphism}, if the lowest
component $F_{0,0}:M_1\to M_2$ induces an isomorphism on homology.
$$ (F_{0,0})_*: H_\bullet(M_1,D^{M_1}_{0,0})\stackrel{\cong}{\to}
H_\bullet(M_2, D^{M_2}_{0,0}).$$
\end{defn}

\begin{thm} \label{invertability}
Let $(BC,D)$ be an A$_\infty$ coalgebra and $F:(B^M
C,D^{M})\to(B^N C,D^N)$ be so that the lowest component $F_{0,0}:
M \to N$ induces an isomorphism $ (F_{0,0})_*: H_\bullet(M,
D^{M}_{ 0,0})\stackrel{\cong}{\to} H_\bullet(N, D^N_{0,0})$
between the homologies. Then, there exists an A$_\infty$
cobimodule morphism $G:(B^N C,D^N)\to (B^M C,D^{M})$ satisfying
$G_*\circ F_*=id_{H_\bullet(B^M C,D^{M})}$ and $F_*\circ
G_*=id_{H_\bullet(B^N C,D^N)}$.
\end{thm}

\begin{proof} For $i=1, 2$, decompose $B^{M_i} C= B^{P_i} C\oplus B^{N_i} C$ into a direct sum of a
minimal and a linear contractible A$_\infty$ cobimodules. Consider
the projection $pr^i:\left(B^{M_i} C,D^{M_i}\right)\to
\left(B^{P_i} C,D^{P_i} \right) $ and the inclusion
$incl^i:\left(B^{P_i} C,D^{P_i} \right) \to \left(B^{M_i}
C,D^{M_i}\right)$. Since $H_*\left( B^{N_i} C,D^{N_i} \right)
=\{0\}$, the maps $pr^i$ and $incl^i$ induce isomorphisms between
the homologies of $(B^{M_i} C,D^{M_i})$ and $(B^{P_i} C,D^{P_i})$.
Let $\Phi=pr^2\circ F \circ incl^1.$
$$
\begin{diagram}
  \node{B^{M_1} C}\arrow{e,t}{F} \node{B^{M_2} C}\arrow{s,r} {pr^2} \\
  \node{B^{P_1} C}\arrow{e,b}{\Phi} \arrow{n,l}{incl^1} \node{B^{P_2} C}
\end{diagram}
$$

Note that $\Phi_{0,0}:(P_1,0)\to(P_2,0)$ induces an isomorphism,
after passing to the homology. Therefore, since the complexes
$(P_1,0)$ and $(P_2,0)$ have zero differentials, $\Phi_{0,0}$ is
in fact an isomorphism. We claim that $\Phi$ has a right inverse
$\rho$, which can be constructed inductively to satisfy $\Phi\circ
\rho=id$. Define $\rho_{0,0}=(\Phi_{0,0}) ^{-1}$. Note that the
$(k,l)$-th level equation $\Phi_{0,0} \circ\rho_{k,l}+
(\text{lower terms in} ~\rho_{r,s})=id$ implies that $\rho_{k,l}=
(\Phi_{0,0})^{-1}\circ (id- \cdots)$. Therefore, $\rho_{k,l}$'s
can be inductively solved. Thus, a right inverse $\rho$ is
constructed. Similarly, one can construct the left inverse
$\lambda$. Thus, $\Phi$ is has the inverse $\rho=\lambda$. Let
$G=incl^1\circ \rho\circ pr^2$. By passing to the homology, we
have,
\begin{eqnarray*}
id_{H_*(B^{M_1}C,D^{M_1})}&=&(incl^1)_*\circ (pr^1)_*\\%
 &=&(incl^1)_*\circ \rho_*\circ\Phi_*\circ(pr^1)_*\\%
 &=&(incl^1)_*\circ \rho_*\circ(pr^2)_*\circ F_*\circ(incl^1)_*\circ(pr^1)_*\\%
 &=& G_*\circ F_*.
\end{eqnarray*}
Similarly, we have,
\begin{eqnarray*}
id_{H_*(B^{M_2}C,D^{M_2})} &=& F_*\circ G_*.
\end{eqnarray*}
\end{proof}

\subsection{Complex of A$_\infty$ co-inner products}

We will now introduce a complex suitable for dealing with
A$_\infty$ co-inner products. By construction, an A$_\infty$
co-inner product will be a closed element in this complex.

Given graded modules $V$, $W$, $X$, and $Y$ over a ring $R$, let
$T^X _Y V=\prod_{k,l\geq 0}V^{ \otimes k} \otimes X \otimes
V^{\otimes l}\otimes Y$. The components are written as $T^X _Y
V_{k,l}=V^{ \otimes k} \otimes X \otimes V^{\otimes l} \otimes Y$
and $T^X _Y V_n=\bigoplus_{k+l=n}V^{ \otimes k} \otimes X \otimes
V^{\otimes l}\otimes Y$. For any map $\mathcal{X}:W\to T^X_Y V$,
denote the component of $\mathcal{X}$ mapping to $T^X_Y V_n$ by
$\mathcal X_n$; that is, $$ \mathcal{X}_n=pr_{T^X_Y V_n}\circ
\mathcal{X}: W \stackrel{\mathcal{X}} {\to} T^X_Y V
\stackrel{pr}{\to} T^X_Y V_n.$$

Let $(BC,D)$ be an A$_\infty$ coalgebra. Recall that A$_\infty$
co-inner products are determined by maps $C^*\to B^C C$. Clearly,
every element $F\in T^C_C C$ gives such a map by dualizing the
last tensor factor. We work with the grading and differential
$\mathcal{D}$ on $T^C_C C$ which makes the inclusion $incl: (T^C_C
C,\mathcal{D})\hookrightarrow (Hom(C^*,B^C C),\delta^{C^*,C})$ is
a chain map of degree zero. To be consistent with the previous
notation, we use the notion $B^C_C C$ for the space $T^C_C C$ with
this correct grading. Therefore, the inclusion $incl: (B^C_C
C,\mathcal{D}) \hookrightarrow (Hom(C^*,B^C C), \delta^{C^*,C})$
is a chain map of degree $0$ with $\delta^{C^*,C}(Im(incl))\subset
Im(incl)$.

Notice that A$_\infty$ co-inner products are exactly the closed
elements of the complex $(Hom(C^*,B^C C), \delta^{C^*,C})$ of
degree zero. Thus, every closed element of $(B^C_C C,\mathcal{D})$
of degree zero determines an A$_\infty$ co-inner product. More
explicitly, the degree of an element,
$$F=(c_1\otimes\cdots \otimes c_n)\otimes c'\otimes
(c''_1\otimes\cdots \otimes c''_m)\otimes c'''\in (C)^{\otimes
n}\otimes C\otimes (C)^{\otimes m}\otimes C,$$ with $c_i\in
C_{k_i}$, $c'\in C_{k'}$, $c''_i\in C_{k''_i}$, and $c'''\in
C_{k'''}$, is given by,
$$ \|F\|=\sum_{i=1}^n(k_i -1)+(k'-1)+
\sum_{i=1}^m(k''_i -1)+(1-k''').$$

The differential $\mathcal{D}:B^C_C C\to B^C_C C$ is given by a
sum obtained by applying $D$ at all possible places in $B^C_C C$
in a cyclic way (see \cite{TT}),
\begin{multline*}
\mathcal D (c_1,\cdots,c_n,c_{n+1},c_{n+2},\cdots,c_{n+m+1},c_{n+m+2})\\
 = \sum_{1\leq i \leq n} \pm (\cdots,c_{i-1},D(c_i) ,c_{i+1},\cdots) \pm
(\cdots,c_{n},D^C (c_{n+1}),c_{n+2},\cdots)\\ + \sum_{n+2\leq i
\leq n+m+1} \pm (\cdots,c_{i-1},D(c_i) ,c_{i+1},\cdots)\pm
\sigma(\cdots, c_{n+m+1},D^C (c_{n+m+2})).
\end{multline*}
where $\sigma:C^{\otimes r}\otimes C\otimes C^{\otimes s}\otimes
C\otimes C^{\otimes t}\to C^{\otimes t+r}\otimes C\otimes
C^{\otimes s}\otimes C=B^C_C C_{t+r,s}$ \label{sig} is a cyclic
rotation of elements. As usual, the signs are determined by the
Koszul sign rule, which says whenever an element of degree $p$
moves over an element of degree $q$, a sign of $(-1)^{pq}$ is
introduced. Let us use a diagrammatic picture for $\mathcal{D}$
described in \cite{TT}. If we draw the two special $C$ components
of $B^C_C C$ on the horizontal axis, then the differential can be
pictured in the following way,
$$ \pm \,\,\,
     \pstree[treemode=R, levelsep=1cm, treesep=0.3cm]{\Tp}
    { \pstree[levelsep=0cm]{
        \Tr
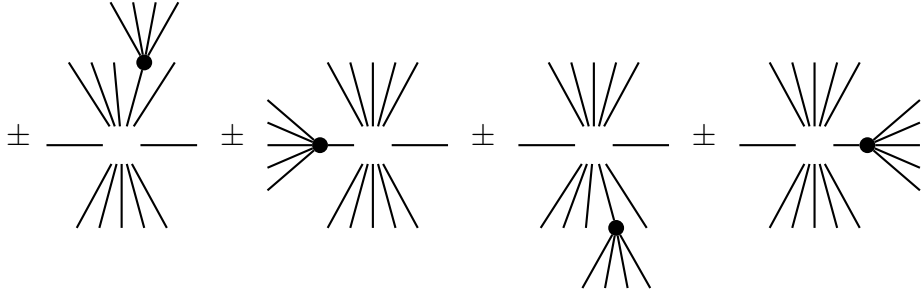
{\begin{pspicture}(0,0)(.5,.5)\end{pspicture} }}
  { \pstree[treemode=U, levelsep=0.8cm]{\Tn}
     {\Tp \Tp \Tp \pstree{\Tc*{3pt}}{\Tp \Tp \Tp \Tp} \Tp}
    \pstree[treemode=R, levelsep=1cm]{\Tn}
     {\Tp}
    \pstree[treemode=D, levelsep=0.8cm]{\Tn}
     {\Tp \Tp \Tp \Tp \Tp}
  }}
 \,\,\, \pm \,\,\,
     \pstree[treemode=L, levelsep=1cm, treesep=0.3cm]{\Tp}
    { \pstree[levelsep=0cm]{
        \Tr{\begin{pspicture}(0,0)(.5,.5)\end{pspicture} }}
  { \pstree[treemode=U, levelsep=0.8cm]{\Tn}
     {\Tp \Tp \Tp \Tp \Tp}
    \pstree[treemode=L, levelsep=0.7cm]{\Tn}
     {\pstree{\Tc*{3pt}}{\Tp \Tp \Tp \Tp \Tp}}
    \pstree[treemode=D, levelsep=0.8cm]{\Tn}
     {\Tp \Tp \Tp \Tp \Tp}
  }}
 \,\,\, \pm \,\,\,
     \pstree[treemode=R, levelsep=1cm, treesep=0.3cm]{\Tp}
    { \pstree[levelsep=0cm]{
        \Tr{\begin{pspicture}(0,0)(.5,.5)\end{pspicture} }}
  { \pstree[treemode=U, levelsep=0.8cm]{\Tn}
     {\Tp \Tp \Tp \Tp \Tp}
    \pstree[treemode=R, levelsep=1cm]{\Tn}
     {\Tp}
    \pstree[treemode=D, levelsep=0.8cm]{\Tn}
     {\Tp \Tp \Tp \pstree{\Tc*{3pt}}{\Tp \Tp \Tp \Tp} \Tp}
  }}
 \,\,\, \pm \,\,\,
     \pstree[treemode=R, levelsep=1cm, treesep=0.3cm]{\Tp}
    { \pstree[levelsep=0cm]{
        \Tr{\begin{pspicture}(0,0)(.5,.5)\end{pspicture} }}
  { \pstree[treemode=U, levelsep=0.8cm]{\Tn}
     {\Tp \Tp \Tp \Tp \Tp}
    \pstree[treemode=R, levelsep=0.7cm]{\Tn}
     {\pstree{\Tc*{3pt}}{\Tp \Tp \Tp \Tp \Tp}}
    \pstree[treemode=D, levelsep=0.8cm]{\Tn}
     {\Tp \Tp \Tp \Tp \Tp}
  }}
$$
(The only difference between the above and \cite{TT} is that here we
apply all of $D=D_1+D_2+\cdots$ including the differential $D_1$,
and not just the higher terms.) Let us recall, why $\mathcal{D
}^2=0$; see \cite{TT}. The diagrams for $\mathcal{D}^2$ are given by
applying $D$ at two places. There are two cases. Either one of the
two multiplications is placed on top of the other, in which case
$D^2=0$ shows that the sum of those diagrams vanish. Or, the
multiplications are placed at different positions in which case each
term appears twice with opposite signs and therefore cancel out each
other.

Given a map $f:BC\to BC'$ between two A$_\infty$ coalgebras
$(BC,D)$ and $(BC',D')$, one can define a map $\hat{f}:B^C_C C\to
B^{C'}_{C'} C'$ by taking a sum over all possibilities of applying
$f$ at all possible places simultaneously. More precisely,
$\hat{f}$ on $C^{\otimes n}\otimes C\otimes C^{\otimes m}\otimes
C$ is given by,
\begin{multline*}
\hat{f}(c_1,\cdots,c_n,c_{n+1},c_{n+2},\cdots,c_{n+m+1},c_{n+m+2})=\\
 \sum_{r,s,t,u} \pm \sigma(f(c_1),\cdots,f(c_n),i^{r,s}(f(c_{n+1})),
 f(c_{n+2}),\cdots,f(c_{n+m+1}),i^{t,u}(f(c_{n+m+2}))).
\end{multline*}
where $i^{r,s}:C^{\otimes r+s+1}\stackrel\cong\to C^{\otimes r}
\otimes C\otimes C^{\otimes s}$ is the canonical isomorphism and
$\sigma$ is the cyclic rotation defined on page \pageref{sig}. In
other words, we have to take the same cyclic rules for the
positions of the elements $c_i$, that were taken in the definition
of $\mathcal {D}$. After applying $f$ in all spots simultaneously,
we need to determine the two special $C$ components. This is done
by taking a sum of all possibilities of special components, which
come from the two original special components,
$$c_{n+1}\mapsto \sum_{k\geq 1} \sum_{r+s+1=k}i^{r,s}
(f_k(c_{n+1})),$$ and,
$$c_{n+m+2}\mapsto \sum_{k\geq 1} \sum_{t+u+1=k}i^{t,u}
(f_k(c_{n+m+2})).$$ In order to be an element of $B^{C'}_{C'} C'$,
the last factor in the tensor product has to be one of the special
$C$ components. Therefore, it might be necessary to apply a cyclic
rotation $\sigma$. Diagrammatically, we have, $$
   \pstree[treemode=D, levelsep=0cm, treesep=0.3cm]
        {\Tr{\begin{pspicture}(0,0)(1.8,0.8)\end{pspicture}}}
  { \pstree[treemode=L, levelsep=0cm]{\Tn}
     {\pstree[treemode=L, levelsep=0.8cm]{\Tr*{f}}{\Tp \Tp \Tp \Tp \Tp}}
    \pstree[treemode=D, levelsep=1cm]{\Tn}{\pstree{\Tr*{f}}{\Tp \Tp \Tp}}
   \pstree[treemode=U, levelsep=1cm]{\Tn}{\pstree{\Tr*{f}}{\Tp \Tp \Tp \Tp}}
    \pstree[treemode=D, levelsep=1cm]{\Tn}{\pstree{\Tr*{f}}{\Tp \Tp \Tp}}
    \pstree[treemode=U, levelsep=1cm]{\Tn}{\pstree{\Tr*{f}}{\Tp \Tp \Tp}}
   \pstree[treemode=D, levelsep=1cm]{\Tn}{\pstree{\Tr*{f}}{\Tp \Tp \Tp \Tp}}
    \pstree[treemode=R, levelsep=0cm]{\Tn}
     {\pstree[treemode=R, levelsep=0.8cm]{\Tr*{f}}{\Tp \Tp \Tp \Tp \Tp} }
  }
$$

\begin{prop}$\hat{f}:(B^C_C C,\mathcal D)\to (B^{C'}_{C'} C',\mathcal D
')$ is a chain map.
\end{prop}

\begin{proof} $\mathcal D '\circ \hat{f}$ corresponds to applying
$f$ and having exactly one multiplication $D'$ outside the ring of
$f$'s; see above picture. $\hat{f} \circ \mathcal D$ corresponds
to applying $f$ and having exactly one multiplication $D$ inside
the ring of $f$'s. But the fact that $f$ is an A$_\infty$ algebra
map ($D'\circ f=f\circ D$) means exactly that $f$ commutes with
inside and outside multiplication of $D'$ and $D$.
\end{proof}

Let us recall the concept of symmetry from \cite{T2} and
\cite{LT}, and its application to the Hochschild complex as stated
in \cite{T2}.

Let $\tau:B^C_C C\to B^C_C C$ denote the map which rotates the
tensor factors cyclically by $180^\circ$; that is,
\begin{equation}\label{tau}
\tau:C^{\otimes n}\otimes C\otimes C^{\otimes m}\otimes C \to
          C^{\otimes m}\otimes C\otimes C^{\otimes n}\otimes C,
\end{equation}
$$ (c_1\otimes\cdots \otimes c_n)\otimes c'\otimes
          (c''_1\otimes\cdots \otimes c''_m)\otimes c'''\mapsto
          (-1)^\epsilon
   (c''_1\otimes\cdots \otimes c''_m)\otimes c'''\otimes
          (c_1\otimes\cdots \otimes c_n)\otimes c',$$
$$
     \pstree[treemode=R, levelsep=1cm, treesep=0.3cm]{\Tr{c' \,}}
    { \pstree[levelsep=0cm]{
        \Tr{\begin{pspicture}(0,0)(.5,.5)\end{pspicture} }}
  { \pstree[treemode=U, levelsep=0.8cm]{\Tn}
     {\Tr{c_n \,\,\,}  \Tr{\cdots} \Tr{\cdots}  \Tr{c_1}}
    \pstree[treemode=R, levelsep=1cm]{\Tn}
     {\Tr{\, c'''}}
    \pstree[treemode=D, levelsep=0.8cm]{\Tn}
     {\Tr{c''_1 \,\,} \Tr{\cdots} \Tr{\cdots} \Tr{\cdots} \Tr{c''_m}}
  }}
  \mapsto (-1)^\epsilon \cdot
     \pstree[treemode=R, levelsep=1cm, treesep=0.3cm]{\Tr{c''' \,}}
    { \pstree[levelsep=0cm]{
        \Tr{\begin{pspicture}(0,0)(.5,.5)\end{pspicture} }}
  { \pstree[treemode=U, levelsep=0.8cm]{\Tn}
     {\Tr{c''_m \,\,} \Tr{\cdots} \Tr{\cdots} \Tr{\cdots} \Tr{c''_1}}
    \pstree[treemode=R, levelsep=1cm]{\Tn}
     {\Tr{\, c'}}
    \pstree[treemode=D, levelsep=0.8cm]{\Tn}
     {\Tr{c_1 \,\,} \Tr{\cdots} \Tr{\cdots}  \Tr{c_n}}
  }},$$ where $\epsilon=(|c_1|+\cdots+|c_n|+|c'|+n+1)\cdot( |c''_1|+
\cdots+
|c''_m|+ |c'''|+m+1)$, and $|c|$ denotes the degree of $c\in C$.

Let $(BC,D)$ be an A$_\infty$ coalgebra. An A$_\infty$ co-inner
product $F\in B^C_C C\subset Hom(B^{C^*}C,B^C C)$ is said to be
\emph{symmetric}, if $F$ is invariant under $\tau$; that is,
\begin{equation}\label{def_symm}\tau(F)=F.
\end{equation}

\begin{defn} \label{duality}
An A$_\infty$ coalgebra with $\infty$ duality consists of an
A$_\infty$ coalgebra $(BC, D)$ together with A$_\infty$ co-inner
product $$F:(B^{C^*}C,D^{C^*}) \to (B^{C}C,D^C),$$ which is
symmetric and induces an isomorphism
$$F_*:H_\bullet(B^{C^*}C,D^{C^*})\stackrel\cong \to
H_\bullet(B^{C}C,D^C),$$ at the level of the homology.
\end{defn}

Let $(BC,D)$ be an A$_\infty$ coalgebra endowed with an $\infty$
duality structure $F$ and strict counit $1\in C$. Denote by
$A=Hom(C,R)$ the dual space of $C$ and let
$CH^\bullet(A,M)=\prod_{n\geq 0}Hom(A^{ \otimes n},M)$ be the
Hochschild cochains of the A$_\infty$ algebra $A$ with values in
an A$_\infty$ bimodule $M$. The following theorem is from
\cite{T2}.

\begin{thm}\label{BV}
Let $F: B^{C^*}C \to B^C C$ denote an $\infty$ duality structure
so that $F_{0,0}$ maps $C^k$ into $C_{d-k}$. Then, $F$ induces an
isomorphism between $HH^\bullet(C^*,C^*)$ and
$HH^{\bullet-d}(C^*,C)$. Moreover, transporting the cup product
from $HH^\bullet(C^*,C^*)$ onto $HH^\bullet(C^*,C)$ and putting it
together with the Connes $\Delta$ operator yields a BV algebra
structure.
\end{thm}

\section{Topological constructions}

\subsection{Construction of an $\infty$ duality}

We now show the relevance of the above algebraic concepts to
topology. Let $X$ be a triangulated space in which the closure of
every simplex is contractible. Let $C=C_\bullet X$ denote the
complex of simplicial chains on $X$. We show that any closed element
$\mu \in C_\bullet X$ gives rise to a symmetric A$_\infty$ co-inner
product $F$. Moreover, if $X$ is a Poincar\'e duality space and
$\mu=[X]$ is its fundamental class, then $F$ is in fact an $\infty$
duality structure. We will then use this to show the Hochschild
cohomology $HH^{\bullet+d}(C^\bullet X, C_\bullet X)$ is a BV
algebra. We construct $F$ in a manner similar to that described by
Dennis Sullivan in Appendix \ref{local-sonstruction-DS}, where a
C$_\infty$ structure on the rational simplicial chains is
constructed.

\begin{defn}
Let $\alpha:C\to C^{\otimes i}$ be a chain map. The map $\alpha$ is called \emph{local}, if it
maps a simplex $\sigma\in C$ into $C(\overline{\sigma})^{\otimes
i}\subset C^{\otimes i}$, where $C(\overline{\sigma}) \subset C$
is the subcomplex of $C$ generated by the cells in the closure of $\sigma$.
\end{defn}

\begin{prop}\label{chain-map-to-inf.i.p.} Let $(BC,D= \widetilde{D_1}+
\widetilde{D_2}+ \cdots)$ be an $A_\infty$ coalgebra on $C$, whose
differential $D$ has local components $D_i:C\to C^{\otimes i}$.
Then, there exists a chain map $\mathcal{X}:(C,D_1) \to(B^C_C
C,\mathcal D)$ of degree $0$, whose lowest component
$\mathcal{X}_0$ is given by $D_2:C\to C\otimes C$.
\end{prop}
\begin{proof}
This is a proof by a double induction on the tensor degree, $n$, and the dimension of the skeleton, $r$.
Thus, we will define local maps $\mathcal{X}_j:C \to B^C_C
C_j$ of degree zero, $j=0,1,2,\cdots$, such that,
\begin{equation}\label{n-level-chi}
 \left(\sum_{0\leq i\leq n} \mathcal{X}_i\right)\circ D_1-
 \mathcal D\circ\left(\sum_{0\leq i\leq n} \mathcal{X}_i\right)
 =\epsilon_{n+1} + \text{higher order terms,}
\end{equation}
where $\epsilon_{n+1}:C \to B^C_C C_{n+1}$ is a map vanishing on
0-simplices and the higher order terms are maps $C \to
\bigoplus_{i>n+1} B^C_C C_i$. Locality means that a simplex
$\sigma\in C$ gets mapped to $B^{C(\overline{\sigma})}
_{C(\overline{\sigma})} C( \overline{\sigma})\subset B^C_C C$.
\begin{enumerate}

\item []${\bf n=0:}$ Let $\mathcal{X}_0=D_2$. This map is local
by assumption and satisfies,
\begin{eqnarray*}
 \mathcal{X}_0\circ D_1-\mathcal D\circ\mathcal{X}_0&=&
  D_2\circ D_1\\
 && -\widetilde{D_1}\circ D_2+ (\text{terms in} ~ \mathcal D
    ~\text{higher than } D_1)\circ D_2\\
 &=& \text{higher order terms.}
\end{eqnarray*}
This is because $D_2$ is a chain map (compare the conditions for $D^2=0$).

\item []${\bf n \geq 1:}$ Assume we have constructed local
maps $\mathcal{X}_j$, $j=0,\cdots,n$ satisfying
(\ref{n-level-chi}). We now start the induction on $r$, the dimension of the skeleton. Thus,
we will construct a local $\mathcal{X}^r_{n+1}$, such that,
\begin{multline}\label{n-r-level-chi}
 \left(\sum_{0\leq i\leq n} \mathcal{X}_i+\mathcal{X}^r_{n+1}\right)\circ D_1-
 \mathcal D\circ\left(\sum_{0\leq i\leq n} \mathcal{X}_i+\mathcal{X}^r_{n+1}\right)
 \\=\epsilon^r_{n+1} + \text{higher order terms,}
\end{multline}
where $\epsilon^r_{n+1}:C \to B^C_C C_{n+1}$ is a map vanishing on
simplices of dimension less or equal to $r$. We may let
$\mathcal{X}^0_{n+1}=0$, because by assumption $\epsilon_{n+1}$
vanishes on $0$-simplices.

Now, assume that a local $\mathcal{X}^r_{n+1}$ has been
constructed. Notice that the map $(\cdot \circ D_1)-(\mathcal
D\circ \cdot):Hom(C,B^C_C C)\to Hom(C,B^C_C C)$ satisfies $((\cdot
\circ D_1)-(\mathcal D\circ \cdot))^2=0$ and therefore,
\begin{eqnarray*}
0&=&((\cdot \circ D_1)-(\mathcal D\circ \cdot))^2
    \left(\sum_{0\leq i\leq n} \mathcal{X}_i+\mathcal{X}^r_{n+1}\right)\\
&=& ((\cdot \circ D_1)-(\mathcal D\circ \cdot))\circ
    ((\cdot\circ D_1)-(\mathcal D\circ \cdot))
    \left(\sum_{0\leq i\leq n} \mathcal{X}_i+\mathcal{X}^r_{n+1}\right)\\
&=& ((\cdot\circ D_1)-(\mathcal D\circ \cdot))
    (\epsilon^r_{n+1} + \text{higher order terms})\\
&=& [D_1,\epsilon^r_{n+1}] + \text{higher order terms.}
\end{eqnarray*}
Therefore, $0=[D_1,\epsilon^r_{n+1}]= \widetilde{D_1} \circ
\epsilon^r_{n+1}- \epsilon^r_{n+1}\circ D_1$. Now, pick any
$(r+1)$-simplex $\sigma$. In order to construct
$\mathcal{X}^{r+1}_{n+1} (\sigma)$, first notice, that
$\epsilon^r_{n+1}$ is a local map because, $$(\epsilon^r_{n+1} +
\text{higher order terms})=((\cdot\circ D_1)-(\mathcal D\circ
\cdot)) \left(\sum_{0\leq i\leq n}
\mathcal{X}_i+\mathcal{X}^r_{n+1}\right),$$ which is local by the
locality assumptions for the $D_i$'s and $\mathcal{X}_i$'s, as
well as the fact that the composition of two local maps is also
local. Thus, for the chosen $(r+1)$-simplex $\sigma$, the map
$\epsilon^r_{n+1}$ restricts to a map
$(\epsilon^r_{n+1})^{\overline{\sigma}}: C(\overline{\sigma}) \to
B^{C(\overline{\sigma})} _{C(\overline{\sigma})}
C(\overline{\sigma})$, which vanishes on simplices of dimension
less than $r+1$. Since $(\epsilon^r_{n+1})^{\overline{\sigma}}$ is
a closed element in $Hom\left(C( \overline{\sigma}), B^{C(
\overline{ \sigma})} _{C(\overline{\sigma})} C(
\overline{\sigma})\right)$ and $\overline{\sigma}$ is by
assumption contractible, we have,

$$H_\bullet\left(Hom\left(C( \overline{\sigma}), B^{C( \overline{
\sigma})} _{C(\overline{\sigma})} C( \overline{\sigma})
\right)\right)=Hom\left(H_0( \overline{\sigma}), H_\bullet
\left(B^{C( \overline{ \sigma})} _{C(\overline{\sigma})} C(
\overline{\sigma})\right)\right).$$

Since $(\epsilon^r_{n+1})^{ \overline{\sigma}}$ vanishes on
$0$-simplices, it is not only closed, but in fact also exact.

Thus, $(\epsilon^r_{n+1})^{\overline{\sigma}}$ can be written as
$(\epsilon^r_{n+1})^{\overline{\sigma}}=- [D_1, \rho^{r+1}_{
n+1}]$, where $\rho^{r+1}_{n+1}:C \to B^C_C C _{n+1}$ can be
assumed to be non-vanishing only on the $(r+1)$-simplex $\sigma$.
We set
$\mathcal{X}^{r+1}_{n+1}=\mathcal{X}^{r}_{n+1}+\rho^{r+1}_{n+1}$
on $C(\overline{\sigma})$, because with this, equation
\eqref{n-r-level-chi} will be satisfied, because,
$$ ((\cdot\circ D_1)-(\mathcal D\circ \cdot)) \left(
\bigoplus_{0\leq i\leq n+1} \mathcal{X}_i+\mathcal{X
}^{r+1}_{n+1}\right) \quad\quad\quad\quad\quad\quad
$$
\begin{eqnarray*} \quad\quad\quad
\quad\quad\quad &=& ((\cdot\circ D_1)-(\mathcal D\circ
\cdot)) \left(\bigoplus_{0\leq i\leq n}\mathcal{X}_i+\mathcal{X}^{r}_{n+1}\right)\\
&&+((\cdot \circ D_1)-(\mathcal D\circ \cdot))(\mathcal{X}^{r+1}_{n+1}) \\
&=&(\epsilon^r_{n+1})^{\overline{\sigma}} + \text{higher order
terms} \\ && +[D_1,\mathcal{X}^{r+1}_{n+1}] + \text{higher order
terms} \\ &=& - [D_1, \rho^{r+1}_{ n+1}]+ [D_1, \rho^{r+1}_{ n+1}]
+\text{higher order terms} \\ &=& \text{higher order terms,}
\end{eqnarray*}
where the higher order terms now have components, $$C(\overline{
\sigma}) \to \bigoplus_{i >n+1} B^{C(\overline{ \sigma})}_{C
(\overline{\sigma})} C(\overline{\sigma})_i.$$ By the
construction, it follows that $\mathcal{X}^{r+1}_{n+1}$ coincides
on the boundary of two $(r+1)$-simplices, and therefore gives rise
a local map $\mathcal{X}^{r+1}_{n+1}$.

Finally, notice, that $\epsilon_{n+2}$ vanishes on $0$-simplices,
since the left-hand side of \eqref{n-r-level-chi} vanishes on
$0$-simplices for any $r$. This completes the induction and the
proof of the proposition.
\end{enumerate}
\end{proof}

\begin{prop}\label{IPonManifold}
Let $X$ be a compact and triangulated space in which the closure
of every simplex is contractible. Let $(BC, D)$ be a local
A$_\infty$ coalgebra structure, where $D_1$ and $D_2$ are the
boundary operator and comultiplication on the on $C=C_\bullet X$,
respectively. For every closed cycle $\mu \in C,$ there exists a
symmetric A$_\infty$ co-inner product $F \in B^C_C C$ such that
the lowest component of $F$ is given by,
$\frac{1}{2}\left(D_2(\mu)+\tau(D_2(\mu))\right)\in C\otimes
C=B^C_C C_0$, where $\tau$ is the rotation map from equation
\eqref{tau}.
\end{prop}
\begin{proof}
Using Proposition \ref{chain-map-to-inf.i.p.}, we obtain a chain map
$\mathcal{X}:(C,D_1)\to (B^C_C C,\mathcal D)$. Let us define an
A$_\infty$ co-inner product by setting $F^0= \mathcal{X} (\mu)\in
B^C_C C$. $F^0$ is in fact an A$_\infty$ co-inner product, since it
is closed under $\mathcal{D}$, $\mathcal D (F^0)=\mathcal D
(\mathcal{X}(\mu))= \mathcal{X}(D_1(\mu))=\mathcal{X}(0)=0$. Now,
let $F=\frac{F^0+\tau(F^0)}{2}$. It is clear that $\tau(F)=F$, since
$\tau^2=id$, and that $\mathcal{D}(F) =0$, since
$\tau\mathcal{D}=\pm\mathcal{D}\tau$.
\end{proof}

\begin{thm}\label{non-deg}
Let $X$ be a compact and triangulated Poincar\' duality space in
which the closure of every simples is contractible and let $\mu
\in C_dX$ denote its fundamental class. Then, the resulting $F\in
B^C_C C\subset Hom(B^{C^*}C , B^C C)$ is an $\infty$ duality
structure. That is to say $F$ is a co-inner product which is
symmetric and induces an isomorphism,
$$ F_*:H_\bullet(B^{C^*} C,D^{C^*})\stackrel\cong \to
H_\bullet(B^{C} C,D^{C}),$$ at the level of the homology.
\end{thm}
\begin{proof}
Choose a strictly cocommutative simplicial chain model
$C=C_\bullet X$ for $X$. This can always be achieved by
symmetrization of a given coproduct. Complete this model to a
local A$_\infty$ coalgebra structure on $C$ as described in
Appendix \ref{local-sonstruction-DS} by Dennis Sullivan.

Recall that in Theorem \ref{IPonManifold} the lowest component
$F_{0,0}:C^* \to C$ of the A$_\infty$ co-inner product is given by
capping a cochain with the fundamental cycle $\mu$, because for
$a, b\in C^*$, it is $(a\otimes b)(D_2(\mu)) =(a\cap \mu)(b)$,
where $D_2$ was chosen to be cocommutative. Since $X$ satisfies
Poincar\'e duality, it follows that $F_{0,0}:C^*\to C$ induces an
isomorphism on homology. Thus, Theorem \ref{invertability} implies
the claim.
\end{proof}

\begin{rem} \label{LocInv} Let us say a few words about a related and important issue.
Although our construction of the A$_\infty$ co-inner product gives
rise to a locally defined structure, the quasi-inverse, which we
show to exist using minimal models, is not necessarily local. This
is simply on account of the fact that the lowest level part of the
inverse map is given by a representative for the Thom class of the
diagonal $X \hookrightarrow X \times X$ and that a result of Mc
Crory \cite{McCrory1977} states the existence of such a local Thom
class implies the Poincar\'e duality space is a homology manifold
(In fact, if the dimension of $X$ is greater than $4$, a result of
Galewski and Stern \cite{G-S} implies that $X$ has the homotopy
type of a topological manifold.) So, as Sullivan pointed out to
us, an $\infty$ duality structure with a local quasi-inverse
should give rise to an $X$-controlled Poincar\'e complex $\it{a}$
$\it{la}$ Ranicki as well as to an L-theory orientation
$[X]_{\mathbb L} \in H_n(X, \mathbb L)$ for $X$. Here, $\mathbb L$
denotes the L-theory spectrum. For an account of algebraic
L-theory and how to obtain rational Pontrjagin classes see
\cite{Ranicki}.
\end{rem}

\subsection{BV structure on ${\bf HH^{\bullet +d}(C^\bullet X, C_\bullet X)}$}

In this chapter, $X$ is a compact and triangulated space
Poincar\'e duality space of dimension $d$, and $\mu \in C=C_dX$
represents its fundamental class. We use the $\infty$ duality
structure of Theorem \ref{IPonManifold} to obtain a BV algebra on
the Hochschild complex $HH^{\bullet+d}(C^\bullet X, C_\bullet X)$.
For this, we use a result from \cite{T2}, namely, that the
Hochschild cohomology of an A$_\infty$ coalgebra with $\infty$
duality is a BV algebra with a unit.

\begin{cor}\label{PD-BV}
Let $X$ be a compact and triangulated Poincar\'e duality space.
Then, the shifted Hochschild cohomology $HH^{\bullet+d} (C^\bullet
X, C_\bullet X)$ is a BV structure with a unit.
\end{cor}
\begin{proof}
This follows immediately by applying the $\infty$ duality structure
obtained in Proposition \ref{non-deg} to Theorem \ref{BV}.
\end{proof}

It was proved in \cite{J} that for a simply connected space $X$
the homology of the free loop space $H_\bullet LX$ is identified
with the Hochschild cohomology $HH^\bullet (C^\bullet X, C_\bullet
X)$. One implication is as follows,

\begin{cor} \label{Loop-BV} Let $X$ be a compact, triangulated, and simply connected
Poincar\'e duality space. Then, the shifted homology of the free
loop space, $H_{\bullet+d} LX$, is a BV algebra with a unit.
\end{cor}

This corollary should be viewed in light of the seminal work of
Chas and Sullivan. It was shown in \cite{CS1} that there exists a
natural BV structure on $H_{\bullet+d} LX$, when $X$ is a closed
manifold. It is therefore natural to ask whether, for a simply
connected and closed smooth manifold, the BV structure in
Corollary \ref{Loop-BV} coincides with that defined in \cite{CS1}.
In fact, developing an algebraic model for the Chas-Sullivan
string topology \cite{CS1} was one of the original motivations for
this work. The identification $H_\bullet LX \cong HH^\bullet
(C^\bullet X, C_\bullet X)$, due to Jones \cite{J}, is explicit.
One can see that the $\Delta$ operator used in string topology
coincides with the $\Delta$ operator on $HH^\bullet (C^\bullet X,
C_\bullet X) $, used in Theorem \ref{PD-BV}. One can ask whether
the product presented here also coincides with the Chas-Sullivan
product. Although there is some evidence that this is indeed the
case, the question remains open. Also, there are identification of
the products in the literature (see \cite{CJ}, \cite{Me}, and
\cite{FTV}), but it is not a priori clear that they are the same
as the one coming from the A$_\infty$ Poincar\'e duality
structure.

\appendix

\section{Local construction of $\infty$ structures \\By Dennis Sullivan}
\label{local-sonstruction-DS}

Let $X$ be a cell complex with cells $e_\alpha$ so that the
closure $\overline{e}_\alpha$ (= $e_\alpha$ union faces of
$e_\alpha$) have the homology of a point ($\mathbb
Q$-coefficients). Let $L(X)=L$ denote the free Lie algebra with
generators $e_\alpha$ placed in degrees (dimension $e_\alpha)-1$.
Consider $L$ as a direct sum $L_0 \oplus L_1 \oplus L_2 \oplus
\cdots $, where $L_0$ is the ground field $\mathbb Q$, $L_1$ is
the linear span of the generators $e_\alpha$, $L_2$ is spanned by
brackets of pairs of elements in $L_1$, etc. Consider derivations
$\delta$ of $L$ expanded into components $\delta=\delta_0 +
\delta_1 + \delta_2 + \cdots $ where $\delta_k$ is determined by a
linear mapping $L_1 \to L_k$ of degree $-1$. Here we discuss only
the special case where $\delta_0=0$.

Orient the cells of $X$ and let $\delta_1$ be the boundary operator
(shifted down by $1$) $L_1 \to L_1$. Note $\delta_1$ still has
degree $-1$. Choose a local canonical cellular approximation to the
diagonal mapping (as in the proof below) to obtain $L_1 \to L_2$.
Extend $\delta_1$ and $\delta_2$ to derivations of $L$ and also call
them $\delta_1$ and $\delta_2$.

\begin{theorem}
There is a canonical local inductive construction of $L_1
\overset{\delta_3}\to L_3$, $L_1 \overset{ \delta_4}\to L_4$,
$\cdots$ so the total derivation $\delta =\delta_1 +\delta_2 +
\cdots$ satisfies $\delta \circ \delta =0$.
\end{theorem}

\begin{remark}
By local we mean $\delta e_\alpha$ belongs in the sub Lie algebra
generated by the cells in the closure of $e_\alpha$.
\end{remark}

\begin{proof}
We interpret $\delta \circ \delta =0$ as the equation $[\delta,
\delta]=0$, where $[\cdot, \cdot]$ is the graded commutator.  For
any $\delta$ the Jacobi identity is $[\delta, [\delta, \delta]]=0$.
Suppose $\delta^k=\delta_1 + \cdots + \delta_k$ has been defined so
that $[\delta^k, \delta^k]$ is the first nonzero term in monomial
degree $k+1$. Jacobi implies this error commutes with $\delta_1$;
that is, it is a closed element in the complex $Der(L)$ of
derivations of $L$. If we work in the closure of a cell, the
homology hypothesis implies that $Der(L)$ has homology only in
degrees $0$ and $-1$. Therefore, the error, which lives in degree
$-2$, can be written as a commutator with $\delta_1$. Using the
cells to generate a linear basis of each $L_k$ by bracketing, we
choose this solution to lie in the image of the adjoint of
$\delta_1$ to make it canonical. This canonical solution is
$\delta_{k+1}$ and this completes the induction, since one knows at
the beginning $\delta_1 \circ \delta_1 =0$ and $\delta_2$ is chain
mapping; that is, $[\delta_2, \delta_1]=0$.
\end{proof}

\begin{remark}\quad
\begin{enumerate}
\item The canonical $\delta$ respects any cellular symmetry of the
cell. \item For the zero cell $a$, then $\delta a + 1/2 [a, a]=0$,
if we choose the diagonal approximation to send a vertex $a$ to
$-a \otimes a$.

\item For the one cell $e$ with vertices $a$, $b$, choose an
orientation so that $\delta_1 e =b-a$ and let $ad_e=[e, \cdot]$.
On the one hand, the procedure of this appendix gives a
differential of the form,
$$\delta{e}=\sum_{i=0}^\infty{\alpha_i}(ad_e)^i(a) +
\beta_i(ad_e)^i(b)\>,$$ and on the other hand, we gave the
following specific infinity structure for the interval in
\cite{LS},
$$\delta{e}=(ad_e)b+\sum_{i=0}^\infty{B_i\over{}i!}(ad_e)^i(b-a)\>,$$
where $B_i$ denotes the $i^{th}$ Bernoulli number. These two
formulae are conjecturally the same.

\item The higher term $\delta_4,\delta_5, \cdots$ have the
interpretation of the rest of the higher homotopies in the
structure of an infinity homotopy co-commutative and
co-associative coalgebra structure on the chains of $X$; that is,
a $C_\infty$ structure.

\item The argument used here can be employed in a variety of
contexts to construct local $\infty$-structures. For instance, let
$C_\bullet X$ denote the complex of chains on $X$ and $TC_\bullet X$
the tensor algebra generated by $C_\bullet X$. Similarly, for a cell
$e_\alpha$, $C_\bullet \overline{e}_\alpha$ denotes the chain
complex of the closure of $e_\alpha$ and $TC_\bullet
\overline{e}_\alpha$ its tensor algebra. Working in the closure of a
cell, $TC_\bullet \overline e_\alpha$ has homology $\mathbb Q$ in
degrees $-1, -2, -3, \cdots$ and therefore $Der(TC_\bullet
\overline{e}_\alpha)$ has homology $\mathbb Q$ in degrees $0, -1,
-2, \cdots$. The error is a class in the homology group
$H_\bullet(Der(TC_\bullet
\overline{e}_\alpha))=Hom(H_0(\overline{e}_\alpha),
H_\bullet(TC_\bullet \overline{e}_\alpha))$. Let us call an element
$\epsilon \in Hom (C_\bullet X, TC_\bullet X)$ local, if for any
cell $\overline{e}_\alpha \subset X$, $\epsilon(C_\bullet \overline
e_\alpha) \subset TC_\bullet(\overline{e}_\alpha)$. We perform a
double induction for both the tensor degree $k$ and the cell degree
$n$. Assume $\delta^{k-1}=\delta_1+\cdots+\delta_{k-1}$ is defined
as a local map on the $n$-skeleton, such that $[\delta^{k-1},
\delta^{ k-1}]=0$ up to tensor degree $k$. Furthermore, assume that
$\delta^{k-1}$ is already extended on the $(n-1)$-skeleton to a
local map $\delta^k= \delta_1+ \cdots+ \delta_{ k-1}+ \delta_k$,
such that $[\delta_k, \delta_k ]=0$ up to tensor degree $k+1$. We
will show that $\delta^k$ can be extended to the $n$-skeleton, such
that $\delta^k$ coincides with the previous map on the faces (i.e.
locality), and $[\delta^k,\delta^k]=0$ up to tensor degree $k+1$. To
this end, look at $\epsilon=[\delta^k,\delta^k]$ on an $n$-cell
$e_\alpha$. $\epsilon$ has a lowest tensor degree of $k$, and
vanishes on the faces of $e_\alpha$. Jacobi identity shows, that
$\epsilon$ is closed under the differential given by the commutator
with $\delta_1$. The complex $Hom(C_\bullet \overline e_\alpha
,TC_\bullet \overline e_\alpha))$ has homology $H_\bullet(Hom(
C_\bullet\overline e_\alpha ,TC_\bullet\overline e_\alpha))
=Hom(H_0(\overline e_\alpha),H_\bullet(TC_\bullet \overline
e_\alpha))$. Now, $\epsilon$ is closed and vanishes on $0$-cells, so
it represents the zero homology class; that is, it must be exact.
Let $\epsilon=[\delta_1, \rho]$. Since $\epsilon$ vanishes on all
faces of $e_\alpha$, we may set $\rho$ equal to zero on all
$m$-cells, for all $m\leq n-1$, and still have the relation
$\epsilon=[\delta_1, \rho]$. Define $\delta_k$ to be the old
$\delta_k$ plus $\rho$. Then, $[\delta^k,\delta^k]=0$ up to tensor
degree $k+1$ on all of $e_\alpha$, and $\delta^k$ is local. To start
the induction, choose $\delta_1$ to be the boundary operator and
$\delta_2$ be a local comultiplication, which is both cocommutative
and coassociative on zero cells. Define $\delta$ on a $0$-cell
by $\delta_k=\Big\{\begin{array}{c} \delta_k \text{, for }k=1,2\\
0\text{, for }k\geq3\,\,\,\end{array}$. Thus, one extends, by
induction on both $n$ and $k$, a local map $\delta$ to the all
cells in the $n$-skeleton, such that $[\delta, \delta]=0$.

\end{enumerate}

\end{remark}
\section{The one dimensional case} \label{circ}

In this section we treat the one dimensional case of the circle. To
this end, we calculate the $\mathcal{X}$ map from Proposition
\ref{chain-map-to-inf.i.p.} for a point, the unit interval, and then
for the circle. We will see, that even in this simplest of
cases, the higher homotopies can not be avoided.

\subsection{Compatibility relations for $\mathcal X$}
Let us start by describing the compatibility relation for
extending the $\mathcal X$ map to the $n$-skeleton once it has
been defined on the $(n-1)$-skeleton of a closed and triangulated
topological space $X$. Using Proposition
\ref{chain-map-to-inf.i.p.}, one can construct a map
$\mathcal{X}:(C,D_1)\to(B^C_C C,\mathcal D)$, where $C=C_\bullet
X$ is a simplicial model of $X$ in which the closure of simplices
are contractible. We will write $A$ for the cochain model
$A=C^\bullet X$. The goal is to build an explicit map
$\mathcal{X}$ inductively over the $n$-skeletons of $C$. We will
sometimes ignore the grading and simply write $A$ instead of the
shifted $sA$, keeping in mind that the correct grading is always
given after shifting.

In detail, $\mathcal{X}:C\to B^C_C C$ being a chain map means that
for all $\sigma\in C$, one has $\mathcal{X}(D_1(\sigma))=\mathcal
D (\mathcal{X}( \sigma))\in B^C_C C\subset Hom(T^A_A A, R)$. Since
$\mathcal{X}(\sigma)\in Hom(T^A_A A, R)$ denotes a sequence of
inner products, we will use the following shorthand,
$$\langle\cdots\rangle^{\sigma}=\mathcal{X}(\sigma):T^A_A
A\to R,
$$ written in components as $
\langle\cdots\rangle^{\sigma}_{k,l}\quad : A^{\otimes k}\otimes A
\otimes A^{\otimes l}\otimes A\to R$. (The arguments from $A$ are
being applied into the dots.) Then, $\mathcal{X}(D_1
(\sigma))=\mathcal D (\mathcal{X}( \sigma))$ reads, $$
\langle\cdots\rangle^{D_1( \sigma)}=\langle\mathcal D
^*(\cdots)\rangle^{\sigma}. $$ Below, we will use the
Alexander-Whitney comultiplication $\Delta:C\to C\otimes C$,
$\Delta(v_0,\cdots,v_n)= \sum_{i=1}^{n} (v_0,\cdots,v_i) \otimes
(v_i,\cdots,v_n)$, which makes $C$ into a strictly coassociative
differential graded coalgebra. In the strictly coassociative case,
the differential $\mathcal D$ only has two components, namely the
differential $D_1$ and the comultiplication $D_2$, because with
this $D=\widetilde{D_1}+\widetilde{D_2}$ satisfies $D^2=0$. Thus,
the chain map condition for $\mathcal{X}$ becomes,
\begin{equation} \label{chain-chi}\langle\cdots\rangle^{D_1(\sigma)}=\langle
\widetilde{D_1^*}(\cdots)\rangle^{\sigma}+\langle\widetilde{D_2^*}(\cdots)\rangle^{\sigma}.
\end{equation}
We will use the notation $\delta(a)=D_1^*(a)$ and $a\cdot b=
D_2^*(a,b)$. The symbol $\widetilde{D_1^*}$ denotes the extension
of $D_1^*$ to a derivation on $T^A_A A$, and similarly
$\widetilde{ D_2^*}$ denotes the extension of $D_2^*$ to a
derivation on $T^A_A A$. It is important to note that two special
elements never combine; see \cite{TT}. For example, if we
underline the two special elements, then we have the maps
\begin{eqnarray*}
 \widetilde{D_1^*} (a,\underline{b},c,\underline{d})&=&
   (\delta(a),\underline{b},c,\underline{d})+(-1)^{\|a\|} (a,\underline{
   \delta(b)},c,\underline{d})\\
  &&+(-1)^{\|a\|+\|b\|} (a,\underline{b},\delta(c),
   \underline{d})+(-1)^{\|a\|+\|b\|+\|c\|}
   (a,\underline{b},c,\underline{\delta(d)}) \\
 \widetilde{D_2^*} (a,\underline{b},c,\underline{d})&=&
  (\underline{a\cdot b},c,\underline{d})+(-1)^{\|a\|} (a,\underline{b\cdot
  c},\underline{d})\\
 &&+(-1)^{\|a\|+\|b\|} (a,\underline{b},\underline{c\cdot d})
  +(-1)^{\epsilon} (\underline{b},c,\underline{d\cdot a} )\\
\widetilde{D_2^*}(a,b,\underline{c},\underline{d})&=&
  (a\cdot b,\underline{c},\underline{d})+(-1)^{\|a\|}
  (a,\underline{b\cdot c},\underline{d})\\
 &&\quad\quad\quad\quad\quad\quad\quad\quad\quad\quad
  \,\,\,\,+(-1)^{\epsilon}
  (b,\underline{c},\underline{d\cdot a})\\
\widetilde{D_2^*}(\underline{a},b,c,\underline{d})&=&
  (\underline{a\cdot b},c,\underline{d})+(-1)^{\|a\|}
  (\underline{a},b\cdot c,
  \underline{d})\\
 &&+(-1)^{\|a\|+\|b\|}(\underline{a},b,\underline{c\cdot d})
\end{eqnarray*}
where $\epsilon=\|b\|+\|c\|+\|a\|\cdot(\|b\|+\|c\|+\|d\|)$, and
$\|x\|$ denotes the degree in $sA$.

Equation \eqref{chain-chi} shows that it is possible to define
$\langle\cdots\rangle ^\sigma$ inductively on the skeleton,
because one only needs to know the lower components
$\langle\cdots\rangle ^{D_1(\sigma)}$. Let $^n\mathcal X$ denote
the $\mathcal X$ map for the standard
$n$-simplex $\Delta^n$. We will inductively define
$^n\mathcal{X}:C_\bullet \Delta^n\to B^{C_\bullet
\Delta^n}_{C_\bullet\Delta^n} C_\bullet \Delta^n$ by extending
the map $^{n-1}\mathcal{X}$ of the standard $(n-1)$-simplex
$\Delta^{n-1}$. If one identifies $\Delta^{n-1}$ with a face of
$\Delta^n$, then $^n\mathcal{X}$ is defined on this face of
$\Delta^n$ by requiring the commutativity of the following
diagram, $$\begin{diagram}
  \node{C_\bullet \Delta^{n-1}}
    \arrow{s,l}{\text{incl}}
    \arrow{e,t}{^{n-1}\mathcal{X}}
  \node{B^{C_\bullet \Delta^{n-1}}_{C_\bullet \Delta^{n-1}} C_\bullet \Delta^{n-1}}
    \arrow{s,r}{\text{incl}} \\
  \node{C_\bullet \Delta^n}
    \arrow{e,b}{^n\mathcal{X}}
  \node{B^{C_\bullet \Delta^n}_{C_\bullet \Delta^n} C_\bullet \Delta^n}
\end{diagram}
$$ It follows, that for a given triangulated topological space $X$
with simplicial chain model $C=C_\bullet X$, one can use the
$^n\mathcal{X}$ to define $\mathcal{X}:C\to B^C_C C$. Namely, for
an n-simplex $\sigma$ of $X$, which is identified with the
standard n-simplex $\Delta^n$, one requires the following diagram
to commute, $$
\begin{diagram}
  \node{C_\bullet\Delta^{n}}
    \arrow{s,l}{\text{incl}}
    \arrow{e,t}{^n\mathcal{X}}
  \node{B^{C_\bullet\Delta^{n}}_{C_\bullet\Delta^{n}} C_\bullet\Delta^{n}}
    \arrow{s,r} {\text{incl}} \\
  \node{C}
    \arrow{e,b}{\mathcal{X}}
  \node{B^C_C C}
\end{diagram}
$$

\subsection{The 0-simplex}
In this section, we calculate the map $\mathcal{X}$ for the trivial
case $\Delta^0=\{a\}$. Let $C=C_\bullet\{a\}$ be of the $0$-simplex,
and $sA=sC^\bullet\{a\}$ be the shifted cochain model, generated by
one element $a$ in degree $1$.

\begin{proposition}\label{Delta0}
If the differential graded algebra structure is given by $\delta(a)
=0$ and $a\cdot a=a$, then $\mathcal{X}(a)=a\otimes a \in C\otimes
C$ defines an A$_\infty$ inner product. More precisely, $\mathcal{X}
(a)$ is given by the following inner products,
\begin{eqnarray*}
\langle a,a\rangle^a_{0,0}&=&1\\
\langle a,\cdots,a\rangle^a_{i,j}&=&0, ~\text{for} ~ i+j > 0
\end{eqnarray*}
\end{proposition}
\begin{proof} We show that the inner products satisfy
equation \eqref{chain-chi}. The first two terms
$\langle\cdots\rangle^{ D_1(a)}$ and
$\langle\widetilde{D_1^*}(\cdots)\rangle^{a}$ always vanish,
because $D_1(a)=0\in C$ and $\delta(a)=0 \in A$. We need to show
that $\langle\widetilde{D_2^*} (\cdots)\rangle^{a}=0$. In order to
check this, we examine $\langle\widetilde{ D_2^*}
(\cdots)\rangle^{a}$ by inserting elements from $T^A_A A$. Notice
that the multiplication takes exactly two inputs and
$\widetilde{D_2^*}$ decreases the monomial degree of an element in
$T^A_A A$ by one, and that the inner product
$\mathcal{X}(a)=\langle\cdots\rangle ^a$ is non-zero only in the
component $A\otimes A=T^A_A A_0 =A^{\otimes 0}\otimes A\otimes
A^{\otimes 0}\otimes A$. Therefore, we need to test equation
\eqref{chain-chi} only for elements of $T^A_A A_1
=\left(A^{\otimes 1}\otimes A\otimes A^{\otimes 0}\otimes A\right)
\oplus \left(A^{\otimes 0}\otimes A\otimes A^{\otimes 1}\otimes
A\right)$, $$ \langle\widetilde{D_2^*} (a,a,a)\rangle_{1,0}^{a}=
   \langle a\cdot a,a\rangle^a_{0,0}-\langle a,a\cdot a\rangle^a_{0,0}=1-1=0 $$
$$ \langle\widetilde{D_2^*} (a,a,a)\rangle_{0,1}^{a}=
   \langle a\cdot a,a\rangle^a_{0,0}-\langle a,a\cdot a\rangle^a_{0,0}=1-1=0 $$
\end{proof}

\subsection{The 1-simplex}
Next, we calculate $\mathcal{X}$ for the case of the interval. Let
$\Delta^1=I=[a,b]$ be an interval, oriented from $a$ to $b$.
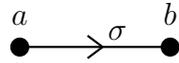
$$
\begin{pspicture}(0,0)(3,1)
 \psline(.5,.5)(2.5,.5)
 \psline(1.4,.65)(1.6,.5) \psline(1.4,.35)(1.6,.5)
 \rput[b](.5,.8){$a$} \rput[b](2.5,.8){$b$}
 \rput[b](1.8,.6){$\sigma$}
 \psdots[dotstyle=*,dotscale=2](.5,.5)
 \psdots[dotstyle=*,dotscale=2](2.5,.5)
\end{pspicture}
$$
Denote by $C=C_\bullet I$ the chain model for the interval given
by the generators $a$, $b$ and $\sigma$. Furthermore, let $sA=s
C^\bullet I$ be the shifted cochain model, which, by slight abuse
of notation, has generators $a$ and $b$ in degree $1$, and
$\sigma$ in degree $0$. The differential graded algebra structure
of $sA$ is given by,
\begin{center}
\begin{tabular}{lll}
$\delta(a)=\sigma,$ & $\delta(b)=-\sigma,$ &$\delta(\sigma)=0,$ \\
\end{tabular}
\end{center}
and,
\begin{center}
\begin{tabular}{lll}
$a\cdot a=a,$ && $b\cdot b=b,$\\
$\sigma\cdot a=-\sigma,$ && $b\cdot \sigma=\sigma,$
\end{tabular}
\end{center}
and all other multiplications vanish. Note that after the shift the
associativity implies $D_2^*(D_2^*(x, y), z)=(-1)^{\|x\|}D_2^*(x,
D_2^*(y, z))$, where $\|x\|$ denotes the degree in $sA$. Notice that
$sA$ is highly non-commutative, e.g., $\sigma\cdot a=-\sigma$, but
$a\cdot\sigma=0$.)
\begin{proposition}\label{Delta1}
Let $\mathcal{X}(a)=a\otimes a$ and $\mathcal{X}(b)=b\otimes b$;
that is,
\begin{eqnarray*}
\langle a,a \rangle^a_{0,0}&=1,\\
\langle b,b \rangle^b_{0,0}&=1,
\end{eqnarray*}
and define $\mathcal{X}(\sigma)$ by the following sequence of inner
products on $sA$,
\begin{center}
\begin{tabular}{ll}
$\langle\sigma,\cdots,\sigma,a\rangle^{\sigma}_{k,0}= 1$& $k\geq 0$\\
$\langle\sigma,\cdots,\sigma,b,\sigma\rangle^{\sigma}_{k,0}=-1$ &
$k\geq 0$
\end{tabular}
\end{center}
and zero otherwise. Then, $\mathcal{X}:C \to B^C_C C$ is a chain
map, where $C=C_\bullet I$.
\end{proposition}
The proof of the above proposition is a case by case calculation
which, in spite of its length, is straightforward and therefore
omitted.

The inner products of $\mathcal{X}(\sigma)$ can be expressed
diagrammatically; see \cite{TT}. For example for $k=5$,
$$
 \begin{pspicture}(0,0.8)(4,2.6)
 \psline(.8,1)(3.2,1)
 \psline(2,1)(1.4,2) \psline(2,1)(0.9,1.6) \psline(2,1)(3.1,1.6)
 \psline(2,1)(2.6,2) \psline(2,1)(2,2.1)
 \psdots[dotstyle=o,dotscale=2](2,1)
 \rput[b](0.7,1.6){$\sigma$}  \rput[b](3.3,1.6){$\sigma$}
 \rput[b](1.4,2.1){$\sigma$}  \rput[b](2.6,2.1){$\sigma$}
 \rput[b](2,2.2){$\sigma$}
 \rput[b](0.5,1){$\sigma$}  \rput[b](3.5,1){$a$}
\end{pspicture} =1
\quad\quad\quad\quad
\begin{pspicture}(0,0.8)(4,2.6)
 \psline(.8,1)(3.2,1)
 \psline(2,1)(1.4,2) \psline(2,1)(0.9,1.6) \psline(2,1)(3.1,1.6)
 \psline(2,1)(2.6,2) \psline(2,1)(2,2.1)
 \psdots[dotstyle=o,dotscale=2](2,1)
 \rput[b](0.7,1.6){$\sigma$}  \rput[b](3.3,1.6){$\sigma$}
 \rput[b](1.4,2.1){$\sigma$}  \rput[b](2.6,2.1){$\sigma$}
 \rput[b](2,2.2){$\sigma$}
 \rput[b](0.5,1){$b$}  \rput[b](3.5,1){$\sigma$}
\end{pspicture} =-1 $$

\begin{remark}
The lowest component of $\mathcal X$ is given by $D_2$. Therefore,
the inner products $\langle b,\sigma\rangle ^{\sigma}$ and
$\langle\sigma,a\rangle^{\sigma}$ are clearly nonzero; see
Proposition \ref{chain-map-to-inf.i.p.}. Let us examine some of
the higher inner products, for example $(b,\underline{b},
\underline{\sigma})$, with special elements underlined.  Equation
\eqref{chain-chi} says,
\begin{eqnarray*}
\langle{b},\underline{b},\underline{\sigma}\rangle^{D_1(\sigma)}_{1,0}&=&\langle
b,\underline{b},\underline{\sigma}\rangle^a_{1,0}-\langle b,\underline{b},\underline{\sigma}
\rangle^b_{1,0}=0 \\
\langle\widetilde{D^*_1}(b,\underline{b},\underline{\sigma})\rangle^{\sigma}_{1,0}&=&\langle\sigma,
\underline{b},
\underline{\sigma}\rangle^{\sigma}_{1,0}-\langle b,\underline{\sigma},\underline{\sigma}
\rangle^{\sigma}_{1,0}\\
\langle\widetilde{D^*_2}(b,\underline{b},\underline{\sigma})\rangle^{\sigma}_{1,0}&=&\langle
\underline{b \cdot b},\underline{\sigma}\rangle^{\sigma}_{0, 0}
=\langle
\underline{b},\underline{\sigma}\rangle^{\sigma}_{0,0} \neq 0 \\
\end{eqnarray*}
This shows that at least one of the two inner products
$\langle\sigma, \underline{b},\underline{\sigma}\rangle^{\sigma}$
or $\langle b,\underline{ \sigma},
\underline{\sigma}\rangle^{\sigma}$ have to be non-zero in order
for equation \eqref{chain-chi} to be satisfied. This example also
shows that there is a choice, which inner products to put
non-zero. We chose $\langle\sigma,
\underline{b},\underline{\sigma}\rangle^{\sigma}\neq 0$, but one
could also choose $\langle b,\underline{ \sigma}, \underline{
\sigma}\rangle^{\sigma}\neq 0$. This choice would result in
another possible A$_\infty$ co-inner product. It turns out that
the choice from Proposition \ref{Delta1} gives the A$_\infty$
co-inner product which has the least amount of non-zero
components.
\end{remark}

After having calculated the map $\mathcal{X}$ for the two lowest
dimensional simplices $\Delta^0$ and $\Delta^1$, it is desirable
to extend it to higher dimensional simplices. Unfortunately, an
explicit calculation of $\mathcal{X}$ for higher simplices becomes
considerably more complicated than the cases considered above.

\subsection{The Circle}
We can now use Proposition \ref{Delta1} to find an A$_\infty$
co-inner product on a model for the circle given by the following
simplicial complex, $$
\begin{pspicture}(0,0)(4,4)
 \pscircle(2,2){1.6}
 \psdots[dotstyle=*,dotscale=2](.4,2)
 \psdots[dotstyle=*,dotscale=2](3.6,2)
 \psline(1.8,3.8)(2,3.6) \psline(1.8,3.4)(2,3.6)
 \psline(2.2,0.6)(2,0.4) \psline(2.2,0.2)(2,0.4)
 \rput[b](0.8,2){$a$} \rput[b](2.2,3.2){$\sigma$}
 \rput[b](3.2,2){$b$} \rput[b](1.8,0.6){$\tau$}
\end{pspicture}
$$
Let $C=C_\bullet S^1$, and $sA=sC^\bullet S^1$ be the shifted
cochain model of the circle with generators $a$, $b$, $\sigma$ and
$\tau$, as indicated above. The differential is given by,
\begin{center}
\begin{tabular}{ll}
$\delta(a)=\sigma - \tau,$ & $\delta(b)=\tau - \sigma,$ \\
$\delta(\sigma)=0,$ & $\delta(\tau)=0,$
\end{tabular}
\end{center}
and the multiplication is non-zero only on the following
generators,
\begin{center}
\begin{tabular}{lcr}
$a\cdot a=a,$ &$\sigma \cdot a = -\sigma,$ &$a\cdot \tau = \tau,$ \\
$b\cdot b=b,$ & $b \cdot \sigma = \sigma,$ & $\tau \cdot b = -\tau.$
\end{tabular}
\end{center}
Again we abuse notation by using the same notation for an element in
$C$ and its dual $A$. The fundamental cycle of $C$ is $\mu= \sigma
+\tau\in C$. Then, by Theorem \ref{IPonManifold}, we know that an
A$_\infty$ co-inner product is given by
$\mathcal{X}(\mu)=\mathcal{X}(\sigma)+\mathcal{X}(\tau)$. Namely,
the only non-zero components of the A$_\infty$ co-inner product are
as follows. For all $k\geq 0$, we have,
\begin{eqnarray*}
\langle\sigma,\cdots,\sigma,a\rangle_{k,0}&=&1,\\
\langle\tau,\cdots,\tau,a,\tau\rangle_{k,0}&=&-1,\\
\langle\sigma,\cdots,\sigma,b,\sigma\rangle_{k,0}&=&-1,\\
\langle\tau,\cdots,\tau,b\rangle_{k,0}&=&1.
\end{eqnarray*}
The following diagrams exhibit the case for $k=5$,
$$\begin{pspicture}(0,0.8)(4,2.6)
 \psline(.8,1)(3.2,1)
 \psline(2,1)(1.4,2) \psline(2,1)(0.9,1.6) \psline(2,1)(3.1,1.6)
 \psline(2,1)(2.6,2) \psline(2,1)(2,2.1)
 \psdots[dotstyle=o,dotscale=2](2,1)
 \rput[b](0.7,1.6){$\sigma$}  \rput[b](3.3,1.6){$\sigma$}
 \rput[b](1.4,2.1){$\sigma$}  \rput[b](2.6,2.1){$\sigma$}
 \rput[b](2,2.2){$\sigma$}
 \rput[b](0.5,1){$\sigma$}  \rput[b](3.5,1){$a$}
\end{pspicture} =1
\quad\quad\quad\quad
\begin{pspicture}(0,0.8)(4,2.6)
 \psline(.8,1)(3.2,1)
 \psline(2,1)(1.4,2) \psline(2,1)(0.9,1.6) \psline(2,1)(3.1,1.6)
 \psline(2,1)(2.6,2) \psline(2,1)(2,2.1)
 \psdots[dotstyle=o,dotscale=2](2,1)
 \rput[b](0.7,1.6){$\tau$}  \rput[b](3.3,1.6){$\tau$}
 \rput[b](1.4,2.1){$\tau$}  \rput[b](2.6,2.1){$\tau$}
 \rput[b](2,2.2){$\tau$}
 \rput[b](0.5,1){$a$}  \rput[b](3.5,1){$\tau$}
\end{pspicture} =-1$$

$$\begin{pspicture}(0,0.8)(4,2.6)
 \psline(.8,1)(3.2,1)
 \psline(2,1)(1.4,2) \psline(2,1)(0.9,1.6) \psline(2,1)(3.1,1.6)
 \psline(2,1)(2.6,2) \psline(2,1)(2,2.1)
 \psdots[dotstyle=o,dotscale=2](2,1)
 \rput[b](0.7,1.6){$\sigma$}  \rput[b](3.3,1.6){$\sigma$}
 \rput[b](1.4,2.1){$\sigma$}  \rput[b](2.6,2.1){$\sigma$}
 \rput[b](2,2.2){$\sigma$}
 \rput[b](0.5,1){$b$}  \rput[b](3.5,1){$\sigma$}
\end{pspicture} =-1
\quad\quad\quad\quad
\begin{pspicture}(0,0.8)(4,2.6)
 \psline(.8,1)(3.2,1)
 \psline(2,1)(1.4,2) \psline(2,1)(0.9,1.6) \psline(2,1)(3.1,1.6)
 \psline(2,1)(2.6,2) \psline(2,1)(2,2.1)
 \psdots[dotstyle=o,dotscale=2](2,1)
 \rput[b](0.7,1.6){$\tau$}  \rput[b](3.3,1.6){$\tau$}
 \rput[b](1.4,2.1){$\tau$}  \rput[b](2.6,2.1){$\tau$}
 \rput[b](2,2.2){$\tau$}
 \rput[b](0.5,1){$\tau$}  \rput[b](3.5,1){$b$}
\end{pspicture} =1$$

\bibliographystyle{amsalpha}

\end{document}